\documentstyle[12pt]{article}

\input{amssym.def}
\input{amssym.tex}

\def\mathfrak{ }

\newcommand{\stp}{{\rm stp}}
\newcommand{\avg}{{\rm Av}}

\newcommand{\qed}{\ \hfill \rule{1ex}{1ex}}

\newcommand{\ci}{\subseteq}

\newcommand{\bprop}{\begin{proposition}}

\newcommand{\eprop}{\end{proposition}}

\newcommand{\bthm}{\begin{theorem}}

\newcommand{\ethm}{\end{theorem}}

\newcommand{\blem}{\begin{lemma}}

\newcommand{\elem}{\end{lemma}}

\newcommand{\bpf}{\begin{proof}}

\newcommand{\epf}{\end{proof}}

\newcommand{\bdefn}{\begin{definition}}

\newcommand{\edefn}{\end{definition}}

\newcommand{\benum}{\begin{enumerate}}

\newcommand{\eenum}{\end{enumerate}}

\newtheorem{definition}{Definition}[section]

\newtheorem{lemma}[definition]{Lemma}

\newtheorem{theorem}[definition]{Theorem}

\newtheorem{proposition}[definition]{Proposition}

\newtheorem{claim}[definition]{Claim}

\newtheorem{fact}[definition]{Fact}

\newenvironment{proof}{{\bf Proof}: }{\hfill \qed \vspace{0.25in}}

\begin{document}

\title{TRANSFERING SATURATION, THE FINITE COVER PROPERTY, AND STABILITY
\thanks{The authors   
thank
the United States - Israel Binational Science foundation 
for supporting this research as well as the Mathematics
department of Rutgers University where part 
of this research was carried out.
}}

\author{John T. Baldwin\thanks{Partially supported by NSF grant 9308768} \\
Department of Mathematics 
\\University of Illinois at Chicago
\\Chicago, IL 60680
\and Rami Grossberg
\\Department of Mathematics
\\Carnegie Mellon University
\\Pittsburgh,  PA 15213
\and Saharon Shelah\thanks{This is item $\#$ 570 on 
Shelah's list of publications.}
\\Institute of Mathematics
\\The Hebrew University of Jerusalem
\\Jerusalem,  91094  Israel\\
\&  \\Department of Mathematics
\\Rutgers University
\\New Brunswick, NJ 08902}\date {\today}
\maketitle
\begin{abstract} 
\underline {Saturation is $(\mu,\kappa)$-transferable in $T$} if and only
if  there is an expansion  $T_1$ of  $T$
with $|T_1| = |T| $
 such that  if
 $M$  is a $\mu$-saturated model of $T_1$  and  $|M| \geq \kappa$
then  the reduct $ M\restriction L(T)$  is $\kappa$-saturated.  We characterize
theories which are
 superstable without f.c.p., or
 without f.c.p. as, respectively those where  saturation is
$(\aleph_0,\lambda)$-transferable or $(\kappa(T),\lambda)$-transferable
for all $\lambda$.  Further if for some $\mu  \geq |T|$, $2^\mu >
\mu^+$, stability is equivalent to for all $\mu \geq |T|$,
saturation is $(\mu,2^\mu)$-transferable.
\end{abstract}


\section{Introduction}

The finite cover property (f.c.p.) is in a peculiar position with respect to
the stability hierarchy.  Theories without the
f.c.p. are stable;  but f.c.p. is independent from $\omega$-stability
or superstability. We introduce a notion of transferability of saturation
which rationalizes this situation somewhat by placing f.c.p. in a natural
hierarchy of properties.
For countable theories the hierarchy is
$\omega$-stable without f.c.p., superstable without f.c.p., not f.c.p.,  
and stable.  
 For appropriate
$(\mu,\kappa)$ each of these
classes of theories is characterized by $(\mu,\kappa)$-transferability
of saturation in following sense.  
\begin{definition}

\underline {Saturation is $(\mu,\kappa)$-transferable in $T$} if and only
if  there is an expansion  $T_1$ of  $T$ with $|T_1|=|T|$  such that  if
$M$  is a $\mu$-saturated model of $T_1$ and $|M| \geq \kappa$,
then the reduct $ M\restriction L(T)$  is $\kappa$-saturated.

\end{definition}

The finite cover property was introduced by Keisler in  \cite{Ke} to 
produce unsaturated ultrapowers. One of his results and a slightly
later set theoretic advance by Kunen yield immediately that if for
$\lambda>2^{|T|}$, saturation is $(|T|^+,\lambda)$-transferable then
$T$ does not have the finite cover property. 
The finite cover property  was  also studied extensively by Shelah in 
\cite{Sh:10} and chapters II, VI and VII of \cite{bible};  those techniques are
used here.

 Our notation generally follows \cite{bible}
with  a few minor exceptions: $|T|$  is the number of symbols in 
$|L(T)|$ plus $\aleph_0$.  We do
not  distinguish between finite sequences and elements,  i.e. we 
write  $a \in A$ to
represent that the elements of the finite sequence  $a$  are from the 
set $A$.  References of the form IV x.y
are to \cite{bible}.


There are several equivalent formulations of the finite cover property.
The following, which looks like a strengthening of the compactness theorem, is most relevant here.

\begin{definition}\label{fcp}
The first order theory $T$ does not have the \underline{finite cover
property} if and only if for every formula $\phi(x;y)$ there exists
an integer $n$ depending on $\phi$ such that for every $A$ contained in
a model of $T$ and every subset $p$ of $\{\phi(x,a), \neg \phi(x,a);
a\in A \}$ the following implication holds:  if every $q \subseteq
p$ with cardinality less than $n$ is consistent then $p$ is consistent.
\end{definition}

The two main tools used in this paper are the following consequence
of not f.c.p. and a sufficient condition for $\lambda$-saturation.

\begin{fact}[II.4.6]\label{easyfact} 
Let  $T$  be a complete first order theory without the f.c.p. and
$\Delta$ a finite set of $L(T)$-formulas.  There is an 
integer $k_{\Delta}$ such that if $M \models T$
is a  saturated
model, $A \ci M$ with $|A| < |M|$ and $ {\bf I}$ is a set of
$\Delta$-indiscernibles over $A$ with cardinality 
at least $k_{\Delta}$
then there exists  ${\bf J} \ci M$   a set of  $\Delta$-indiscernibles (over $A$)
extending $ {\bf I}$ of cardinality  $|M|$.
\end{fact}

The principal tool for establishing the transfer of saturation is

\begin{fact}[III.3.10]\label{basictool} 
If a model $M$ of a stable theory
 if $M$ is either $F^a_{\kappa(T)}$-saturated  or 
$\kappa(T) +\aleph_1$-saturated and
for each set of indiscernibles ${\bf I}$ contained in $M$ there is an
equivalent set of indiscernibles ${\bf J}$ contained in $M$
with $|{\bf J}| = \lambda$ then
$M$ is $\lambda$-saturated.
\end{fact}  

We thank Anand Pillay for raising the issue of the superstable case
and the  referee for 
the final formulation of Theorem~\ref{nfcp} which
generalizes our earlier version and for correcting an oversight
in another argument.

\section{The transferability hierarchy }

In this section we characterize 
certain combinations of stability and the finite cover property
in terms of transferability of saturation. 
Extending the notation we write {\em saturation is
$(0,\kappa)$-transferable in $T$} if and only
if  there is an expansion  $T_1$ of  $T$
with $|T_1|=|T|$  such that  if $M \models T_1$ and $|M| \geq \kappa$,
 $M|L(T)$
is $\kappa$-saturated.  In particular, taking $|M| = \kappa$, 
$PC(T_1,T)$ is categorical in $\kappa$.  Using  this language
we can reformulate an old result of Shelah (
Theorems VI.5.4 and VIII.4.1) to provide the first stage of our hierarchy.

\begin{fact}
\label{fact1}  For a countable theory $T$, the following are equivalent.
\begin{enumerate}
\item $T$ does not have the finite cover property and is $\omega$-stable.
\item For all $ \lambda > \aleph_0 $, saturation is 
$(0,\lambda)$-transferable in $T$.
\item  For some $ \lambda > 2^{\aleph_0}$      
 saturation is $(0,\lambda)$-transferable in $T$.
\end{enumerate}
\end{fact}

Since the proof of 1) implies 2) is not given in \cite{bible}
and follows the line of our other arguments we sketch the proof
in our discussion after Theorem~\ref{ss}.   This result holds
only for countable languages; the remainder 
apply to theories of arbitrary cardinality.

We introduce the following special notation to uniformize the
statement of the next result.

\[\kappa'(T) = \left\{ \begin{array}{ll}
\kappa(T) & \mbox{ if $T$ is stable}\\
|T|^+  &  \mbox{ if $T$ is unstable}
\end{array}
\right.
 \]

\begin{theorem}
\label{thmA}   
The following are equivalent for a complete theory $T$.
\begin{enumerate}
\item $T$ does not have the finite cover property.
\item For all $ \lambda \geq |T|^+ $, saturation is 
$(\kappa'(T),\lambda)$-transferable in $T$.
\item For some $ \lambda > 2^{|T|}$,    
saturation is $(\kappa'(T),\lambda)$-transferable in $T$.
\end{enumerate}
\end{theorem}

\begin{proof}
It is obvious that $(2)$ implies $(3)$.  Now we show $(3)$ implies $(1)$ 
by showing that if $T$ has the  f.c.p. then saturation is not
even $(|T|^+,\lambda)$-transferable 
(and so certainly not $(\kappa(T) +\aleph_1,\lambda)$-transferable). 
Let $T_1$ be any extension of $T$ and
$N_0$ an arbitrary model of $T_1$ with cardinality at least 
$\lambda$.   By Kunen's theorem (see \cite{Ku}, or Theorem 6.1.4 in \cite{CK})
there exists an $\aleph_1$-incomplete $|T|^+$-good ultrafilter $D$  on  $|T|$.
Denote by $N_1$ the ultrapower $N_0^{|T|}/D$.  By \cite{Ke} 1.4  and 4.1 or VI.5.3,
$N_1$ is $|T|^+$-saturated but not  $(2^{|T|})^+$-saturated.

We now show 
 $(1)$ implies $(2)$.
Let  $T$  be a theory without the f.c.p..  By  II.4.1,  $T$  is stable.
The proof now splits into two cases depending on whether $T$ is
superstable.  We begin with the case that $T$ is stable but not superstable.
Then $\kappa'(T) \geq \kappa(T) + \aleph_1$ and this inequality
will be essential shortly.

Let   $L_1:=L(T) \bigcup \{F\}$  where $F$  is a binary function symbol.  
The theory  $T_1$
consists of $T$  and the following axioms.
\begin{enumerate}
\item  For each  $x$,  the function  $F(x,\cdot)$  is injective.
\item  For every finite  $\Delta \ci L(T)$,  let  $k_{\Delta}$  
be the integer from
Fact \ref{easyfact}.  If
$I$ is a finite set of $\Delta$-indiscernibles of cardinality 
at least  $k_{\Delta}$
then there
exists an $x_I$  such that \begin{enumerate}
\item the range of  $F(x_I,\cdot)$ contains  $I$  and
\item the range of  $F(x_I,\cdot)$  is a  set of $\Delta$-indiscernibles.
\end{enumerate}
\end{enumerate}

It should be clear that the above axioms can be formulated in 
first order logic in the
language $L_1$.
To see that $T_1$ is consistent, we expand a saturated model  $N$
of $T$  to a model of $T_1$.
Fix a $1$-$1$ correspondence between finite sets
of $\Delta$-indiscernibles {\bf I} with  $| {\bf I}| \geq k_{\Delta}$
and elements $x_{\bf I}$ of $N$. By
Fact \ref{easyfact}, each sufficiently large 
finite sequence of $\Delta$-indiscernibles {\bf I}
in $N$ extends to one with $|N|$ elements.  
Fix a $1$-$1$ correspondence between the universe of $N$ and this sequence.
Interpret
$F(x_I,x)$ as this correspondence.

Suppose $N^*$ is a $\kappa'(T)$-saturated model of $T_1$
of cardinality at least $\lambda$ and denote
$N^*\restriction L(T)$ by $N$.
 Since $\kappa'(T) \geq \kappa(T) +\aleph_1$,
by Fact~\ref{basictool} we need only establish the following claim.

\begin{claim}
\label{claim2}  Any infinite sequence of indiscernibles
${\bf  I}$  in $N$ extends to a sequence {\bf J} of
indiscernibles (over the empty set)
with cardinality $|N|$. 
\end{claim}

\begin{proof}
Let $ q(x)$ be the  set of formulas
which expresses that
for each finite $\Delta$
 the range of $F(x,\cdot)$
is a set of $\Delta$-indiscernibles and ${\bf I}$ is contained in the
range of $F(x,\cdot)$.
If $ a \in N$
realizes the type  $q(x)$ then, since $F(a,\cdot)$ is $1$-$1$,   
${\bf J}:= \{ F(a,b) : b \in N \}$  is as required.
We now
show $q(x)$  is consistent.   
Fix a finite $q^* \ci q(x)$ and let  $\Delta $  be a
finite subset of $ L(T)$ such that
all the $L(T)$-formulas from $q^*$  appear in $\Delta$. 
Let  $m< \omega$  be sufficiently large
so that all the elements of  $I$ appearing in $q^*$ are among 
$\{b_0, \ldots ,b_{m-1} \}$  and
$m \geq k_{\Delta}$. It suffices to show that for some $a\in N$,
each $b_i$ for $i<m$ is in the range of $F(a,\cdot)$ and the range
of $F(a,\cdot)$ is a set of $\Delta$-indiscernibles.
This follows immediately from  $T_1$, 
by the assumption that $m
\geq k_{\Delta}$. Since $q$ is over a countable set
 there exists an element  $a \in N^*$  satisfying $q^*$ and we finish.
\end{proof}$_{\ref{claim2}}$

We now prove Case 2) of $1$ implies $2$:  superstable $T$.  The
general outline of the proof is the same but we replace $\kappa(T)
+ \aleph_1$-saturation with $F^a_{\kappa(T)}$-saturation and we must use
a different trick to find an equivalent set of  indiscernibles.
The idea for guaranteeing 
$F^a_{\kappa(T)}$-saturation is taken from Proposition 1.6 of \cite{Ca};
the referee suggested moving it from a less useful place in the
argument to here.  

\begin{lemma}
\label{kT}
If $T$ does not have the f.c.p. then there is an expansion $T_1$ of
$T$ such that if $M$ is a $\kappa(T)$-saturated model of $ T_1$ 
then $M\restriction L(T)$ is
$F^a_{\kappa(T)}$-saturated.
\end{lemma}

\begin{proof}
Let  $T$  be a theory without the f.c.p..  
Form $L_1$ by adding to $L$  new $k$-ary function symbols $f^{\theta,E}_i$, 
for $i<m= m(\theta,E)$,
for each pair of   formulas $\theta(z), E(x,y,z)$ with $\lg(z) =k$
 such that
for any $M \models T$ and $a \in M$, if $M\models \theta(a)$ then
$E(x,y,a)$ is an equivalence relation with $m$ classes.
The theory  $T_1$
consists of $T$  and the following axioms:
 For each $k$-ary sequence $z$ such that $\theta(z)$, 
the elements $f^{\theta,E}_i(z)$, $i<m$
provide a complete set of representatives for $E(x,y,z)$.
In any model of $T$, one can choose Skolem functions $f^{\theta,E}_i(z)$ to
give sets of representatives for the finite equivalence relations so $T_1$
is consistent.  Now suppose that $N^* \models T_1\;$ is  $\kappa(T)$-saturated.
For any $q =\stp(d/C)$ with
$|C|<\kappa(T)$, note that $q$ is equivalent to the $L_1$-type 
over $C$ consisting of the formulas
$E(x,f^{\theta,E}_i(c))$ for $E$ a finite equivalence relation defined
over a finite sequence $c \in 
C$ such that $E(d,f^{\theta,E}_i(c))$.  Since this type is realized,
$N=N^*\restriction L(T)$ is $F^a_{\kappa(T)}$-saturated.
\end{proof}$_{\ref{kT}}$

Now we show  finish showing $(1)$ implies $(2)$ in the superstable case. 
Let $\lambda \geq |T|^+$  be given.  
 We must find a $T_2$ to 
witness $(\aleph_0,\lambda)$-transferability.
First expand $T$ to $T_1$ as in Lemma~\ref{kT} so that
 if $M$ is an $\aleph_0$-saturated model of $ T_1$ 
then $M\restriction L(T)$ is
$F^a_{\aleph_0}$-saturated.
Form $L_2$ by adding to $L_1$ an $n+2$-ary function symbol $F_n$ for
each $n$. 
The theory  $T_2$
consists of $T_1$  and the following axioms:
\begin{enumerate}
\item  For each  $x$ and $n$-ary sequence $z$,  the function  
$F_n(x,z,\cdot)$  is injective.
\item  For every finite  $\Delta \ci L(T)$ and $n$-ary sequence $z$,  
let  $k_{\Delta}$  
be the integer from
Fact \ref{easyfact}.  If
$I$ is a finite set of $\Delta$-indiscernibles over $z$ of cardinality 
at least  $k_{\Delta}$
then there
exists an $x_I$  such that 
\begin{enumerate}
\item the range of  $F_n(x_I,z,\cdot)$ contains  $I$,
\item the range of  $F_n(x_I,z,\cdot)$  
is a  set of $\Delta$-indiscernibles
over $z$.
\end{enumerate}
\end{enumerate}

The consistency of $T_2$ is entirely analogous to the similar step
in the proof of 
Theorem~\ref{thmA}.  We just have to interpret each $F_n(x,z,y)$
instead of a single function of two variables.
Now suppose that $N^* \models T_2\;$ is an $\aleph_0$-saturated model  
of cardinality at least $\lambda$. 
The reduct $N$ of $N^*$  to $L(T)$ is
$F^a_{\aleph_0}$-saturated and it suffices by Fact~\ref{basictool}
to show
for each set of indiscernibles $I$ contained in $N$ there is an
equivalent set of indiscernibles $J$ with $|J| = \lambda$.

Let ${\bf I} =\{b_n : n< \omega \}$ be such an infinite set of 
indiscernibles in
$N$.  Let $p^* = \avg(I,N)$ and, since $N$ is $F^a_{\aleph_0}$-saturated,
choose  $m < \omega$ such that for  $B = \{b_0 \ldots b_{m-1}\}$,
  $p^*|B$ is stationary and
$p^*$ does not fork over $B$.
We show 
there is  a sequence {\bf J} of
indiscernibles based on $p^*\restriction B$
with $|{\bf J}| = |N|$.  
Let $q_1(x)$ be a type over $B$ that contains 
$(\forall y)\theta(F_m(x,b_0, \ldots b_{m-1},y))$ for all $\theta(x) \in
p^*\restriction B$,  for each $\phi(x_0,\ldots x_{n-1}) \in L(T)$ such that
$N \models \phi(b_0, \ldots b_{n-1})$, the formula
$ (\forall y_1) \ldots (\forall y_n) 
\phi( F_m(x,b_0, \ldots b_{m-1},y_1),\ldots F_m(x,b_0, \ldots b_{m-1},y_n))$
and  the assertion that $F_m(x,c,\cdot)$ is injective.
The definition of $T_2$ implies the consistency of $q_1$. 
Since $q_1$ is a type over a finite set,  $q_1$ is realized by
 an element  $a \in N^*$;  
this
guarantees the existence of
a set of $|N|$ indiscernibles equivalent to $I$ as required.
\end{proof}$_{\ref{thmA}}$

In the superstable case we can get one slightly stronger result which
allows to characterize superstable without f.c.p. by
$(\aleph_0,\lambda)$-transferability. 

\begin{theorem}\label{ss}   
If for some  $\lambda > 2^{|T|}$,
saturation is $(\aleph_0,\lambda)$-transferable in $T$
 then $T$ is superstable without the f.c.p. 
\end{theorem}

\begin{proof}
By Theorem~\ref{thmA}
 we deduce  that $T$
does not have the f.c.p.~ (using $\lambda > 2^{|T|}$)  
and so
$T$ is stable.  Suppose for contradiction there are a stable but not superstable
$T$ and a $T_1$ which witnesses $(\aleph_0,\lambda)$-transferability in $T$.  
Apply VIII.3.5 to ${\rm PC}(T_1,T)$ taking $\kappa = \aleph_0$, 
$\mu = (2^{|T|})^+$ and $\lambda \geq\mu$.  There are $2^{\mu}$  models
of $T_1 $  with cardinality $\lambda$, which are $\aleph_0$-saturated, whose
reducts to $L(T)$ are nonisomorphic.  So some are not $\lambda$-saturated.
\end{proof}$_{\ref{ss}}$.

 We were unable to find a
uniform argument for $1$ implies $2$ of Theorem~\ref{thmA}; there seem to be two quite
different ideas for making the large set of indiscernibles equivalent
to the given set.  The proof of 1) implies 2) of Fact~\ref{fact1}
proceeds along similar lines with the following variation.  Since $T$ is
$\omega$-stable every $\omega$-saturated model is $F^a_{\omega}$-saturated.
Again using the $\omega$ stability, it easy to Skolemize with countably
many functions so that each type over a finite set is realized.  Then
the same trick as in Theorem~\ref{ss} guarantees the existence of large
equivalent
indiscernible sets.

The proof of Theorem~\ref{thmA} yields somewhat more than is necessary. 
The theory $T_1$ 
which is found in the implication (1) implies (2) does not depend on 
$\lambda$ and contains only a single additional function symbol.
We could obtain a stronger result than $(3)$ implies $(1)$
with the same proof by demanding in a modified $(3)$
that the model witnessing $(|T|^+,\lambda)$-transferability have cardinality
$\lambda = \lambda^{|T|} > 2^{|T|}$.
 
As pointed out by the referee,
we can combine the arguments for Theorem~\ref{thmA} and Theorem~\ref{ss}
to characterize $ \kappa(T)$ for theories without
the finite cover property if  $\kappa(T)$
satisfies the set-theoretic conditions of Theorem VIII.3.5.  For example,
under the GCH  if $\kappa(T)$ is not the successor of 
a singular cardinal and $T$ does not have the f.c.p. $\kappa(T)$ is
the least $\kappa$ such that there is $\lambda > 2^{|T|}$ for which
saturation is $(\kappa,\lambda)$-transferable.

\begin{theorem}\label{stable}
 Suppose that there exists a cardinal  
$\mu \geq |T|$  such that
$\;2^{\mu}>\mu^+$.  
For a complete theory   $T$, the following are equivalent:
\begin{enumerate}
\item  $T$  is stable.
\item For all $ \mu \geq |T|$,
 saturation is $(\mu^+,2^{\mu})$-transferable in $T$.
\item For some $ \mu \geq |T|$,  saturation is
 $(\mu^+,\mu^{++})$-transferable in $T$.
\end{enumerate}
\end{theorem}

The condition $\mu^+ < 2^{\mu}$ is used only for $(2)$ implies $(3)$ 
(which is obvious with that hypothesis).  In the 
next two lemmas we prove 
in ZFC that $(1)$ implies $(2)$,
and that $(3)$ implies $(1)$.  This shows in ZFC that stability 
is bracketed between two transferability conditions.

\begin{lemma}
\label{B1->2}
If $T$  is stable and
$  \mu \geq |T|$, saturation is $(\mu^+,2^{\mu})$-transferable in $T$.
\end{lemma}

\begin{proof} 
We must find an expansion $T_1$ of  $T$  such that if 
$M$ is a $\mu^+$-saturated
model of $T_1$ and $|M| \geq 2^{\mu}$, 
$M\restriction L$ is $2^{\mu}$-saturated.  Form $L_1$ by
adding one additional binary predicate $E(x,y)$ and add axioms asserting 
that $E$ codes all finite sets.  (I.e., for every set of $k$ elements 
$x_i$ there is a unique $y$ such that $E(z,y)$ if and only if $z$ is 
one of the $x_i$.)  For any model $M_1$ of $T_1$ and any element $b$ 
of $M_1$, let $[b]:=\{a\in M_1 : M_1\models E[a,b]\}$.

Now let $M_1$ be a $\mu^+$-saturated model of $T_1$ and $M$ the 
reduct of $M_1$ to $L$.  Suppose $A \subseteq M$ has cardinality less 
than $2^{\mu}$ and $p \in S^1(A)$.  We must show $p$ is realized in 
$M$.  By the definition of $\kappa (T)$  there exists $B\subseteq A$
of cardinality less than $\kappa (T)$  such that $p $ does not fork 
over $B$.  Since $M_1$ is $|T|^+$-saturated, 
we may take $p \restriction B$ to be stationary.
Let ${\hat p}\in S(M)$ be an extension of $p$ that does not fork over 
$B$.
Since $\mu^+>|T|\geq\kappa (T)$, by  $\mu^+$-saturation of $M$ there 
exists
$I:=\{a_n : n<\omega\}\subseteq M$ such that
$a_n\models {\hat p}\restriction (B \cup \{a_k : k<n\})$.  Since
the sequence is chosen over a stationary type, 
$I$ is a set of indiscernibles.

Define an $L_1$-type  $q$ over $I$ so that if $b$ realizes $q$,
$[b] \bigcup I$ is a set of indiscernibles over the empty set.  
Since the relation $E$ codes finite sets, and $I$ is a set of 
indiscernibles
$q$ is consistent.  By the $\aleph_1$-saturation of $M_1$  there 
exists
$b\in M$ realizing
the type $q$. If 
$[b]$ has $2^{\mu}$ elements we are finished since for each formula 
$\phi({ x},{\overline y})$ and each
${\overline a} \in A$ with $\phi({ x},{\overline a})\in p$, only 
finitely many elements of $[b]$
satisfy $\neg\phi({ x},{\overline a})$.
To show $[b]$ is big enough, using the $\mu^+$-saturation of $M$, 
we define inductively
for $\eta \in 2^{\leq\mu}$ elements $c_{\eta}\in M$ such that
\begin{enumerate}
\item $c_{\emptyset}=b$
\item For any $\eta$, $[c_{\eta\frown 0}]$ and $ [c_{\eta\frown 1}]$
are disjoint subsets of $[c_{\eta}]$.
\item If $\lg(\eta)$ is a limit ordinal $\alpha$, 
$[c_{\eta}] \ci \cap_{i<\alpha}[c_{\eta\restriction i}]$
\end{enumerate}

Now for $s \in 2^{\mu}$, the $c_s$ witness that $|[b]| = 2^{\mu}$.
\end{proof}$_{\ref{B1->2}}$.

\begin{lemma}
\label{B3->1}
If $ \mu \geq |T|$ and  
saturation is $(\mu^+,\mu^{++})$-transferable in $T$
then $T$ is stable.
\end{lemma}

\begin{proof}
 For the 
sake of contradiction suppose $T$ is an unstable theory and 
that 
there is a $T_1 \supseteq T$ such that 
if $M$ is a $\mu^+$-saturated model of  $T_1$
with cardinality at least $\mu^{++}$,
$M\restriction L(T)$ is $\mu^{++}$-saturated.
Fix $M_0\models T_1$ with cardinality
at least $\mu^{++}$.  Let $D$ be a $\mu$-regular ultrafilter
on $I=\mu$.   Construct an ultralimit sequence
$\langle M_{\alpha}:\alpha < \mu^+\rangle$
as in VI.6
with $M_{\alpha+1} = M_{\alpha}^I/D$ and taking unions
at limits.  By VI.6.1  $M_{\mu^+}$
is $\mu^+$-saturated.  But by
VI.6.2, since $T$ is unstable, $M_{\mu^+}$ is not
$\mu^{++}$-saturated.
\end{proof}$_{\ref{B3->1}}$


The methods and concerns of this paper are similar to those
in the recent paper of E. Casanovas \cite{Ca}.
He 
defines a model to be expandable if every consistent expansion of 
${\rm Th}(M)$ with at most  $|M|$ additional symbols can be realized 
as
an expansion of $M$.  His results are orthogonal to those here.  He 
shows for countable stable $T$ that  $T$ has an expandable model 
which is not saturated of cardinality greater than the continuum if 
and only if $T$ is not superstable or $T$ has the finite cover 
property.

By varying the parameters in $(\mu,\kappa)$-transferability of saturation
we have characterized four classes of countable theories:
$\omega$-stable without f.c.p., superstable without f.c.p., not f.c.p., and
stable.  For uncountable $\lambda$, they 
correspond respectively to:
$(0,\lambda)$-transferability, $(\aleph_0,\lambda)$-transferability,
$(\aleph_1,\lambda)$-transferability,
$(\aleph_1,2^{\aleph_0})$-transferability.  Although the analogous
results for uncountable languages are more cumbersome to summarise, 
countability of the language
is only essential for the $\omega$-stable characterization.


\begin{thebibliography}{99}
\bibitem[Ca]{Ca} Enrique Casanovas, Compactly expandable models and stability, {\em Journal of Symbolic
Logic} {\bf 60} , 1995, pages  673--683.
\bibitem[CK]{CK}   C.C. Chang and H. Jerome Keisler,  {\bf Model Theory}, 
North-Holland Pub.l Co.  1990.
\bibitem[Ke]{Ke} H. Jerome Keisler,  Ultraproducts which are not saturated, {\em Journal of Symbolic
Logic} {\bf 32} , 1967, pages  23--46.
\bibitem[Ku]{Ku} Kenneth Kunen,  Ultrafilters and independent sets,  {\em Trans. Amer. Math. Soc.},
{\bf 172}, 1972, pages 199--206.
\bibitem[Sh 10]{Sh:10}Saharon Shelah,
\newblock {Stability, the f.c.p., and superstability; model theoretic
properties of formulas in first order theory},
\newblock {\em {Annals of Mathematical Logic}}, {\bf 3}, pages 271--362, 1971.
---  {\bf MR:}~47:6475, (02H05)
\bibitem[Sh c]{bible}  Saharon Shelah, {\bf Classification Theory and 
the Number
of Nonisomorphic Models}, Rev. Ed., North-Holland, 1990, Amsterdam.
\end{thebibliography}
\end{document}